\documentclass{amsart}
\title[Multizeta values, period interpretation and relations]{Multizeta values 
for $\FF_q[t]$, their period interpretation and relations between them}
\author[Anderson]{Greg W. Anderson}
\address{University of Minnesota, Mpls., MN 55455}
\email{gwanders@umn.edu}
\author[Thakur]{Dinesh S.\ Thakur}
\address{University of Arizona, Tucson, AZ 85721}
\email{thakur@math.arizona.edu}
\date{January 10, 2009}
\newcommand{\ZZ}{{\mathbb{Z}}}
\newcommand{\CC}{{\mathbb{C}}}
\newcommand{\OO}{{\mathcal{O}}}
\newcommand{\NN}{{\mathbb{Z}_{\geq 0}}}
\newcommand{\FF}{{\mathbb{F}}}
\newcommand{\pp}{\tilde{\pi}}

\newtheorem{Theorem}{Theorem}

\begin{document}

\begin{abstract}
We provide a period interpretation for 
 multizeta values (in the
function field context) in terms of explicit iterated extensions of
tensor powers of Carlitz motives (mixed Carlitz-Tate $t$-motives). We give
examples of combinatorially involved relations that these multizeta
values satisfy.

\end{abstract}

\maketitle

\setcounter{section}{-1}
\section{Introduction}

The multizeta values introduced and studied originally by Euler have
been pursued again recently with renewed interest because of their
emergence in studies in mathematics and mathematical physics
connecting diverse viewpoints. 
They occur naturally as coefficients of the Drinfeld associator, and
thus have connections to quantum groups, knot invariants and
mathematical physics. They also occur in the Grothendieck-Ihara
program to study the absolute Galois group through the fundamental
group of the projective line minus three points and related studies of
iterated extensions of Tate motives, Feynman path integral
renormalizations, etc. We refer the reader to papers on this subject
by Broadhurst, Cartier, Deligne, Drinfeld, \'{E}calle, Furusho,
Goncharov, Hoffman, Kreimer, Racinet, Terasoma, Waldschmidt, Zagier,
Zudilin to mention just a few names.

Having learned about these rich interconnections at the Arizona Winter
school, the second author, in 2002, defined and studied two types of
multizeta values 
\cite[Sec 5.10]{T} for function fields, one complex-valued
(generalizing special values of Artin-Weil zeta functions) and the other with values in
Laurent series over finite fields (generalizing Carlitz zeta values).
For the $\FF_q[t]$ case, the first type was completely evaluated in
\cite{T} (see \cite{M} for 
a study in the higher genus
case),
for both types some identities were established,
and for the second type failure of the shuffle identities was noted.  Because of the failure,
some other variants of the second type were also investigated in \cite{T}.
Here we deal  only with the second type, and we restrict attention 
exclusively to $\FF_q[t]$.

We now outline the main points  of the paper.
(i) We introduce `degenerate multizeta values', thus bringing
back the failed sum shuffle relations, and in turn making the span of the
degenerate and original (here called non-degenerate) multizetas into an
algebra. (ii) We construct explicit 
mixed Carlitz-Tate $t$-motives with periods evaluating to
multizeta values (generalizing the results for zeta values
obtained in \cite{AT}), and we construct enough of them to account for all
multizeta values.
(iii) We give some examples of new combinatorially involved multizeta
identities, discuss their $t$-motivic nature, and 
indicate contrasts with the classical case. (For a more systematic account
of such identities, see \cite{T2}.)
Finally, (iv) without proving anything, 
we discuss the transcendence problems to which our work can be applied.
It is our hope 
that, just as the construction given in \cite{AT} helped 
determine \cite{CY} all
the relations among the Carlitz zeta values, 
the constructions given in this paper will help to
determine  all relations among multizeta values.  This hope 
is a key motivation for the paper.

 Recent ongoing work \cite{T2} 
suggests that the introduction of degenerate multizeta values is not
really necessary in the sense that the original multizeta value span is also an
algebra, on account of 
new kinds of combinatorially involved identities. We
hope that eventually the latter will be understood as 
analogs of integral shuffles once the relevant theory is developed.

Finally, let us 
state our main motivation:
to provide a few clues toward the discovery
of the $t$-motivic analog of the fundamental group
of the projective line minus three points. Let's call it 
provisionally $t$-$\pi_1$. We hope to discover or inspire
others to discover it through detailed investigation of its
various children. The mixed Carlitz-Tate $t$-motives constructed here are
presumably pro-nilpotent 
children of $t$-$\pi_1$.
In other work \cite{AT2} we have found what plausibly could
be 
metabelian children of $t$-$\pi_1$,
 i.e., reasonable analogs of Deligne-Soul\'{e} cocycles
 and Ihara power series, and we have linked these objects with the special
  points constructed in \cite{AT} and \cite{Aprime}. 
  We hope that
 the $t$-$\pi_1$-heuristic leads next to the analog of the Drinfeld associator.

\section{Multizeta values over $\FF_q[t]$}
\subsection{Notation}
$$\begin{array}{rcl}

\ZZ&=&\mbox{\{integers\}}\\
\ZZ_+&=&\mbox{\{positive integers\}}\\
 \NN&=&\mbox{\{non-negative integers\}}\\
T,\tilde{t}&=&\mbox{independent variables}\\
q&=&\mbox{a power of a prime} \ p\\
t&=&-\tilde{t}^{q-1}\\
A&=&\FF_q[t]\\
A_+&=&\mbox{monics in $A$}\\
K&=&\FF_q(t)\\
\tilde{K}&=&\FF_q(\tilde{t})\\
K_\infty&=&\FF_q((1/t))=\mbox{completion of $K$ at $\infty$}\\
\tilde{K}_\infty&=&\FF_q((1/\tilde{t}))=\mbox{completion of $\tilde{K}$ at $\infty$}\\
\OO_\infty&=&\FF_q[[1/\tilde{t}]]\subset \tilde{K}_\infty\\
\CC_\infty&=&\mbox{completion of algebraic closure of $\tilde{K}_\infty$}\\
\mbox{[$n$]}&=&t^{q^n}-t\\
D_n&=&\prod_{i=0}^{n-1}(t^{q^n}-t^{q^i})\\
\ell_n&=&\prod_{i=1}^n(t-t^{q^i})\ = \ (-1)^nL_n \\
\deg&=&\mbox{function assigning to $x\in K_{\infty}$ its degree in $t$}
\end{array}
$$

\subsection{Carlitz  zeta values and power sums}

\subsubsection{}
For $s\in \ZZ_+$, put
$$\zeta(s)=\sum_{a\in A_+}\frac{1}{a^s}\in K_\infty.$$
These are the {\em Carlitz zeta values}. 
See \cite{G, T} and references there for more. 

\subsubsection{}
It is convenient to break the Carlitz zeta values into power
sums grouped by degree: 
Given integers $s>0$ and $d\geq 0$ put
$$S_d(s)=
\sum_{\begin{subarray}{c}
a\in A_+\\
\deg a=d
\end{subarray}}\frac{1}{a^s}.$$
(This is $S_d(-s)$ in the notation of \cite{T}).

\subsection{Multizeta values and multiple power sums}

\subsubsection{Multizeta values}\label{subsubsection:Official}
For $s=(s_1,\dots,s_r)\in \ZZ_+^r$, following \cite[Sec. 5.10]{T}, we define the multizeta value 
$$\zeta(s)=\zeta(s_1, \cdots,
s_r)=\sum_{d_1>\cdots >d_r\geq 0}
S_{d_1}(s_1)\cdots S_{d_r}(s_r)=\sum
\frac{1}{a_1^{s_1}\cdots
  a_r^{s_r}}\in K_\infty,$$ where the second sum is over $(a_1,\dots,a_r)\in A_+^r$ 
  with $\deg a_i$ strictly decreasing.  We say that this multizeta value has {\em depth} $r$ and {\em weight} 
$\sum s_i$.
\subsubsection{Remark}
In \cite[Sec. 5.10]{T} the notation
  $\zeta_d(s)$ was used, where the subscript $d$ was just meant
to call the word ``degree'' back to mind.

\subsubsection{Multiple power sums}
As with the Carlitz zeta function, it is convenient to break the sum defining 
a multizeta value down according to degrees. To that end we introduce the
following notation.
Given $r$-tuples $d=(d_1,\dots,d_r)\in \NN^r$
and $s=(s_1,\dots,s_r)\in \ZZ_+^r$, put
$$S_d(s)=S_d(s_1, \cdots, s_r)=\prod_{i=1}^r S_{d_i}(s_i).$$
Then we have
$$\zeta(s)=\sum S_d(s),$$
where the sum runs through all $d=(d_1,\dots,d_r)\in (\NN)^r$ with 
the $d_i$'s strictly
decreasing.

\subsubsection{Motivation for degenerate multizeta values}
\label{subsubsection:MotivationDegenerate}

Now classically, sum shuffle identities prove that the multizeta
values 
span an algebra over $\ZZ$.  The algebra structure is something we
obviously want to preserve on the function field side.  But as shown
in \cite{T} (see also Section 3 below), in the function field
situation, you cannot
usually shuffle $a$'s in $A_+$.  Ultimately the
problem is that whereas $\ZZ$ has a natural total order, $A$ has none.
But all is not lost.  We may choose to shuffle degrees rather than
individual elements of $A_+$.  Then, at the expense of somewhat
broadening the definition of multizeta values, we once again get a
full set of sum shuffle identities, enough to prove that the
$\FF_p$-linear combinations of $\FF_q[t]$-multizeta values form an
algebra over $\FF_p$ (see \ref{subsubsection:FFSumShuffle} below for
the explicit identity proving this).

These considerations bring us to the following definition generalizing that
given in \ref{subsubsection:Official}.

\subsubsection{Degenerate multizeta values}\label{subsubsection:OfficialDegenerate}
For each subset $I$ of $\{1, 2, \cdots, r-1\}$ and \linebreak $s\in \ZZ_+^r$, we define 
$$\zeta_I(s):= \sum S_d(s)\in K_{\infty}$$
where we sum over $d=(d_1, \cdots, d_r)\in (\NN)^r$ with monotonically
decreasing $d_i$'s having jumps exactly in positions in $I$. 
Of course for $I=\{1,\dots,r-1\}$ we have $\zeta(s)=\zeta_I(s)$.
For $I\neq \{1,\dots,r-1\}$ we call $\zeta_I(s)$ a {\em degenerate}
multizeta value. For contrast we say that $\zeta(s)$ is 
a {\em nondegenerate} multizeta value.

\subsubsection{Remarks}
(i) See 3.6 for an explanation of why we may not really need 
these degenerate multizeta values to get an algebra.

(ii) Alternatively, we could consider  a variant of degenerate 
multizeta values defined by conditions of the type $d_1>d_2\geq d_3 >\cdots$
in place of conditions of the type $d_1>d_2=d_3>\cdots$ as above. These 
are clearly simple $\FF_p$-linear combinations of the 
degenerate multizeta values already defined.

\section{Period interpretation}

\subsection{Twisting}
Given 
$$f=\sum_{i=0}^\infty f_iT^i\in \CC_\infty[[T]],$$
put
$$f^{(n)}=\sum_i f_i^{q^n}T^i\in \CC_\infty[[T]].$$
Extend this rule entry-wise to matrices with entries in $\CC_\infty[[T]]$.

\subsection{The function $\Omega$}
$$\Omega(T)=\tilde{t}^{-q}\prod_{i=1}^\infty (1-T/t^{q^i})
\in \OO_\infty[[T]].$$
This has an infinite radius of convergence, and satisfies the
functional equation
$$\Omega^{(-1)}(T)=(T-t)\Omega(T).$$ 
The quantity 
$$\pp=1/\Omega(t)$$ is a period
of the Carlitz module (\cite[p.\ 179]{AT} or \cite{T}) and could be
thought of as the analog of $2 \pi i$. ($\pp$ is exactly the same  notation as in \cite{T}.)

\subsection{Carlitz gamma and factorial} 
Given $n\in \NN$, we define the {\em Carlitz gamma} $\Gamma$ and {\em factorial} $\Pi$ by
$$\Gamma_{n+1}:=\Pi_n:=\prod_i D_i^{n_i},$$
where
$$n=\sum_{i=0}^\infty n_iq^i\;\;\;(0\leq n_i\leq q)$$
is the base $q$ expansion of $n$.

\subsection{The polynomials $H_s$ and  the period representation of Carlitz zeta}

\subsubsection{Interpolation of power sums}
 From \cite[3.7.4]{AT},  we know there exists
for each $s\in \NN$ a unique polynomial $H_s=H_s(T)\in A[T]$
such that
\begin{equation}\label{equation:ZetaInterpolationFormula}
(H_{s-1}\Omega^s)^{(d)}|_{T=t} = \Gamma_s S_d(s)/\pp^s
\end{equation}
for all $d\in \NN$ and $s\in \ZZ_+$.
It is known that $H_s\Omega^{s+1}\in \pi \OO_\infty[[T]]$
for all $s\in \NN$. 

\subsubsection{} Explicitly, $H_s(T)$ is defined by the generating series 
identity

\begin{equation}\label{equation:ProShtuka}
\sum_{s=0}^{\infty}\frac{H_s(T)}{\Gamma_{s+1}|_{t=T}}x^{s}=
\left(1-\sum_{i=0}^{\infty}
\prod_{j=1}^{i} \frac{(T^{q^i}-t^{q^j})}{(T^{q^i}-T^{q^{j-1}})} 
x^{q^i}\right)^{-1}. 
\end{equation}

\subsubsection{Remarks}
(i) Our $H_s(T)$ is the same as the two variable polynomial $H_s(y,
T)$ of \cite{AT} evaluated at $y=t$.  (ii) Our notation here neatly
(or confusingly, depending on
your point of view) interpolates between 
\cite{ABP} and \cite{T} in the sense that the latter two references
use $t$ and $T$ in exactly opposite ways. And note that, instead of
$T$, the letter $\theta$ was used in \cite{A}.  Unfortunately the
reader of this literature just has to live with a jumble of $t$'s,
$T$'s and the occasional $\theta$.

\subsubsection{Delayed interpolation of power sums}\label{subsubsection:Delayed} Combining 
the interpolation formula above with the functional equation in 2.2, 
we see that, for $w\geq 1$, 
$$(H_{s-1}^{(-w)}[(T-t)^{(0)}\cdots (T-t)^{(-w+1)}]^s\Omega^s)^{(d+w)}|_{T=t}
=\Gamma_s S_d(s)/\pp^s.
$$

\subsubsection{} 
In \cite[3.8.2]{AT}, using the interpolation formula 
\eqref{equation:ZetaInterpolationFormula}, we constructed  an 
algebraic point on the $s$-th tensor power
 $C^{\otimes s}$ of the Carlitz module $C$ whose logarithm connected
 to the Carlitz zeta value $\zeta(s)$. 
 This is equivalent to constructing an
 extension over $A$ of  $C^{\otimes n}$ by the trivial module $C^{\otimes 0}$ having
 $\Gamma_s\zeta(s)$ as its period.

\subsection{Iterated extensions and multizeta values as periods}
Let us first see the mechanism of how iterated sums as in the
definition of multizeta sums come up naturally as entries in triangular
matrices satisfying shtuka $F-1$ equations. 
The original
$t$-motive formalism \cite{A} or the equivalent dual
$t$-motive formalism used in \cite{ABP} will then realize these
matrices as period matrices of appropriate iterated extensions of
tensor powers of the Carlitz module. These should be viewed as analogs of
iterated extensions of Tate motives $\ZZ(s)$.

\subsubsection{} Let $s=(s_1,\dots,s_r)\in \ZZ_+^r,$. Let $[[X]]$
denote the diagonal matrix $$[X^{s_1+\cdots +s_r}, X^{s_2+\cdots +s_r},
\cdots, X^{s_r}, X^0]. $$
Put $\Lambda=[[\Omega]]$ and $D=[[(T-t)]]$ and 
$$Q=\left[\begin{array}{ccccc}
1\\
Q_{21}&1\\
&\ddots&\ddots\\
&&Q_{r+1,r}&1
\end{array}\right],\;\;\;
L=\left[\begin{array}{cccc}
1\\
L_{21}&\ddots\\
\vdots&\ddots&\ddots\\
L_{r+1,1}&\dots&L_{r+1,r}&1
\end{array}\right],
$$

$$\Phi=Q^{(-1)}D,\;\;\;\Psi=\Lambda L.$$
For the moment we leave $L_{ij}$ and $Q_{ij}$ undefined. 
These we will determine
presently.

\subsubsection{} We have then 
$\Lambda^{(-1)}=D\Lambda$ and so the uniformizability relation (see below) 
\begin{equation}\label{equation:UniformizabilityRelation}
\Phi \Psi=\Psi^{(-1)}, 
\end{equation}
is equivalent to
\begin{equation}\label{equation:UniformizabilityRelationBis} 
(\Lambda^{-1}Q\Lambda)L^{(1)}=L.
\end{equation}

\subsubsection{}
We have
$$\Lambda^{-1}Q\Lambda=\left[\begin{array}{ccccc}
1\\
\Omega^{s_1}Q_{21}&1\\
&\ddots&\ddots\\
&&\Omega^{s_r}Q_{r+1,r}&1
\end{array}\right], 
$$ 
and hence the relation \eqref{equation:UniformizabilityRelationBis}  translates into recursions 
$$L_{ij}=\Omega^{s_{i-1}}Q_{i,i-1}L_{i-1,j}^{(1)}+L_{ij}^{(1)}, 
\ \ i>j+1$$
$$L_{i,i-1}=\Omega^{s_{i-1}}Q_{i,i-1}+L_{i,i-1}^{(1)}.$$

\subsubsection{}  Let us now solve the recursions. We have
$$
\begin{array}{cc}
L_{21}-L_{21}^{(1)}=\Omega^{s_1}Q_{21},&
L_{21}=\sum_{i=0}^\infty (\Omega^{s_1}Q_{21})^{(i)},\\
\vdots&\vdots\\
L_{r+1,1}-L_{r+1,1}^{(1)}=\Omega^{s_r}Q_{r+1,r}L_{r,1}^{(1)},&
L_{r+1,1}=\sum_{i=0}^\infty (\Omega^{s_r}Q_{r+1,r}L_{r,1}^{(1)})^{(i)},
\end{array}
$$
with entries on the right  telescoping to entries on the left,
and then
$$L_{r+1,1}= \sum_{ i_1>\dots>i_r\geq 0}(\Omega^{s_r}Q_{r+1,r})^{(i_r)}
\cdots (\Omega^{s_1}Q_{21})^{(i_1)},$$
as a consequence of the the lines on the right.

We thus see (after evaluation at $t=T$) the iterated sum expression of 
the type defining multizeta
values.
This calculation due to the first author was the starting point of this 
paper. 

By the theory of $t$-motives \cite{A, G, T} and the equivalent theory
of dual $t$-motives \cite{ABP, P}, the period matrix for the dual
$t$-motive defined by $\Phi$, with $\Psi$ related to it as in the
uniformizability relation \eqref{equation:UniformizabilityRelation} is
given by $\Psi^{-1}$ evaluated at $T=t$.  (For the reader familiar
with $t$-motives but not dual $t$-motives, we just note here that in
passing from $t$-motives to dual $t$-motives, the residues in the
recipe for $t$-motive periods in \cite{A}, \cite[7.4]{T} are changed
to evaluation.)

\subsubsection{}
Since $\Psi=\Lambda L$, we have $\Psi_{ij}=\Omega^{s_i+\cdots
 +s_r}L_{ij}$.  
\subsubsection{}
If we let 
\begin{equation*}\label{equation:QijDef}
Q_{i+1, i}=H_{s_i-1}(T)
\end{equation*}
 then we see that $L_{r+1, 1}$ evaluated at $T=t$ is
$$\Gamma_{s_1}\cdots \Gamma_{s_r}\zeta(s_1, \cdots, s_r)/\pp^{s_1+\cdots
+s_r}.$$ 

\subsubsection{} In fact, for $i>k\geq 1$, we see that $\Psi_{ik}$
evaluated at $T=t$ is $$\Gamma_{s_k}\cdots
\Gamma_{s_{i-1}}\zeta(s_k, \cdots, s_{i-1})/\pp^{s_k+\cdots
+s_{r}}.$$

\subsubsection{} Let us give a simpler recipe for the 
entries of $\Psi=(\psi_{ij})$ 
and the period matrix 
$$\Psi^{-1}=(p_{ij}).$$

Let us write temporarily $Z_{ij\cdots}$ for $\Gamma_{s_i}\Gamma_{s_j}\cdots
\zeta(s_i, s_j, \cdots)$. Given an expression $F$ in $Z_{ij\cdots}$, by 
$F(k)$ we denote the corresponding expression obtained from $F$ by increasing 
all the indices by $k$. For example, $F=Z_1Z_2-Z_{12}$ would
give $F(1)=Z_2Z_3-Z_{23}$.

The $j$-th column of $\Psi$ is $\pp^{-(s_j+\cdots +s_r)}$ times 
$[0,\cdots, 0, 1, Z_j, Z_{j(j+1)}, \cdots , Z_{j(j+1)\cdots r}]$. 

Let us `normalize'  as follows:
$$\psi_{ij}'=\psi_{ij}/\pp^{-(s_j+\cdots+s_r)},\;\;\;
p_{ij}'=p_{ij}/\pp^{s_i+\cdots+s_r}.$$

Then $p_{(i+k)(j+k)}' =p_{ij}(k)' $ and $\psi_{(i+k)(j+k)}'
=\psi_{ij}'(k)$. So entries of these normalized matrices 
are immediately obtained from those of the first column 
by shifting and adding indices. We have already described 
the first column for (normalized) $\Psi$. Let us describe 
it for the (normalized) period matrix by giving its first column. 
\begin{eqnarray*}
p_{i1}'&=&-(p_{i-1,1}'(1)Z_1+p_{i-2,1}'(2)Z_{12}+\cdots +
p_{11}'(i-1)Z_{1,2,\dots,i-1})\\
&=&-(p_{i2}'Z_1+\cdots+p_{ii}'Z_{1,2,\dots,i-1}).
\end{eqnarray*}
In particular, note that the bottom-left entry is $-Z_{12\cdots r}$ 
plus (or minus) products of lower depth $Z$ entries (without any 
$\pp$-powers) and $\prod \Gamma_{s_i}$ is a common factor. 
 For example, in rank 2, we get $Z_1$. In  rank 3, we get
$Z_1Z_2-Z_{12}=\Gamma_{s_1}\Gamma_{s_2}[\zeta(s_1)\zeta(s_2)-\zeta(s_1, s_2)]$.
In rank 4, we get $Z_1Z_{23}+Z_{12}Z_3-Z_1Z_2Z_3-Z_{123}$.

\subsubsection{} Since $H_s^{(-1)}(T)\in A[T]$, as follows from its
explicit description above, the entries of $\Phi$ are defined over
$A[T]$.  Thus the non-degenerate multizeta value $\zeta(s)$ (times
$\prod \Gamma_{s_i}$) occurs as an entry of the inverse of the period
matrix of mixed Carlitz-Tate $t$-motive over $A$.

Let us summarize  what we have proved: 

\begin{Theorem}
  With $Q_{i, j}$ defined as in 2.5.6, the (dual) $t$-motive over $A$
  defined in 2.5.1 is a rank $r+1$ motive with the inverse of the
  period matrix given as in 2.5.7, in particular containing the depth $r$
  multizeta value $\zeta(s_1,\dots,s_r)$.
\end{Theorem}

\subsubsection{} 
Recall the degenerate multizeta values defined in 
\ref{subsubsection:OfficialDegenerate}.
To take care of these values in a slight expansion of the framework of the preceding theorem,
we modify as follows.
We use the polynomials in the delayed interpolation formula of \ref{subsubsection:Delayed} in place 
of the polynomials $H_{s-1}$ whenever there are no jumps of degrees.
In this way we get  $t$-motives which have arbitrary degenerate multizeta values
(multiplied by appropriate gamma factors) as periods.
However, the matrices $\Phi$ arising this way have entries  in $\FF_q[t^{1/q^{w}}][T]$.
In other words, we need to make an inseparable base extension to realize those values.

\subsubsection{Remarks} (i) Bloch and Terasoma inform us that in the
classical case, it is not known whether multizeta values of depth $r$
can be achieved as periods of rank $r+1$ motives, although given a
multizeta value there is a combinatorial way using the description of
Grothendieck-Teichm\"uller group to find the rank.  In the function
field situation we can say at least that depth $r$ multizeta values
appear in the {\em inverse} of the period matrix of a dual $t$-motive
of rank $r+1$.

(ii) From the calculations above, we see that multizeta value of depth
$r$ plus a linear combination of lower depth multizeta values appears
as an entry of the period matrix.  To make each multizeta value of depth
$r$ itself appear as an entry of the period matrix (rather than an
entry of the inverse), we need only form tensor products and direct
sums of the basic dual $t$-motives already constructed.  For example,
(ignoring gamma factors by working over $K$ rather than $A$), since
$\zeta(s_1)\zeta(s_2)-\zeta(s_1, s_2)$ occurs in rank 3, and
$\zeta(s)$ in rank 2, we can get $\zeta(s_1, s_2)$ in rank 7 (taking
direct sum of rank 3 with tensor product of the two rank 2).

(iii) For some classical treatments in the number field case, the results and
the period calculations exist only over $\mathbb{Q}$ and not over $\ZZ$.
In contrast, our approach in the function field situation naturally gives 
results over $\FF_q[t]$.

(iv) It was shown by Abhyankar that in finite characteristic any curve
(not necessarily defined over an algebraic extension of the prime
field) can be realized as an \'{e}tale cover of the affine line, and
so one has in principle a Belyi-type embedding of the absolute Galois
group of $\FF_q(t)$ into the outer automorphism group of the algebraic
fundamental group of the affine line over $\overline{\FF_q(t)}$. But
on account of wild ramification the latter group is overwhelmingly
complicated!  Accordingly, we have bypassed the Belyi-embedding
approach to deal directly with the relevant mixed motives. (Indeed, it
is not clear that $t$-$\pi_1$ need be a fundamental group.)  In the
depth one case, a better justification for our {\em ad hoc} approach
is obtained through construction and study of close analogs \cite{AT2}
of the Deligne-Soul\'{e} cocycles, the Ihara power series and related
structures.  (In fact, the power series \eqref{equation:ProShtuka} in
a suitable sense defines a ``pro-t-motive'' giving rise to the Ihara
power series analog.)  The higher depth situation is a work in
progress and will be dealt with elsewhere.

\section{Relations between multizeta values} 

First,
until 3.6, we focus on some basic facts about non-degenerate
multizeta values and refer to \cite{T2} for more general results.

\subsection{Relations to the $p^{th}$-power map}
Since we are in characteristic $p$, the definition immediately implies that 
the multizeta value at $(ps_i)$ is the $p$-th power of the corresponding 
multizeta value at $(s_i)$. 

\subsection{Zeta at `even' $s$}
Carlitz proved \cite[5.2.1]{T} an  analog of Euler's result 
that if $s$ is `even' in the sense that $q-1$ divides 
$s$, then $\zeta(s)/\pp^s$ is in $K$. 

\subsection{Low $s$ relations}
Using $S_d(kp^n)=1/\ell_d^{kp^n}$, for $0<k<q$ \cite[5.9.1]{T}, it was
noted in \cite[5.10.6]{T} that the non-degenerate multizeta values with
weight not more than $q$ satisfy classical sum-shuffle identities.
More generally, $\prod_j \zeta(s_i(j))$ is a sum of multizeta values as in
classical sum shuffle identities, if all sums over $j$ of any
$s_i(j)$'s is not more than $q$. 

So any classical sum-shuffle relation with fixed 
$s_i(j)$'s works for $q$ large enough. 

Note here though that because of the characteristic, the look and
consequences of these relations can be different. We have, by 3.3,
$\zeta(1)\zeta(1)=\zeta(2) +2\zeta(1, 1)$ working for all $\FF_q[t]$
as $q$ varies, but the the right side reduces to $\zeta(2)$ when $p=2$ and
the relation does not give any information on $\zeta(1, 1)$. Also,
when $p=2$, but not in general, the sum shuffle
$\zeta(k)\zeta(k)=\zeta(2k)+2\zeta(k, k)$ works by 3.1, but the right side
reduces to $\zeta(2k)$ and the relation does not give any information
on $\zeta(k, k)$.

\subsubsection{} More generally, if  $s_i=ap^m$, $s_{i+1}=bp^n$, 
$s_i+s_{i+1}=cp^k$, with $a, b, c\leq q$, then $S_d(s_i)S_d(s_{i+1})
=S_d(s_i+s_{i+1})$, and thus a degenerate multizeta value, with 
no jump at an index $i$,  equals the multizeta value of 
one depth lower obtained by putting $s_i+s_{i+1}$ in combined $i^{th}$ and 
$(i+1)^{st}$ place.

\subsection{Failure of naive sum shuffle }

In ~\cite[Thm. 5.10.12]{T}, it is shown that when $q=3$, 
$\zeta(2, 2)/\pp^4$ is not 
in $K$. So  when $q=3$, using  Carlitz's result above for 
`even' $s$, we see that the naive analog of 
classical sum shuffle $\zeta(2)^2=2\zeta(2, 2)+\zeta(4)$ fails. 
(See below for a much simpler proof of the last fact).

\subsection{Failure of naive integral  shuffle }

Classically, there are integral shuffle identities between multizeta values 
coming from their iterated integral expressions and thus connecting 
immediately 
to mixed motives. For example,  in the usual iterated integral notation we have 
$$\int w_1w_0\int w_1w_0 =2\int w_1w_0w_1w_0 +4 \int w_1w_1w_0w_0$$
which with $w_0=dz/z$ and $w_1=dz/(1-z)$ gives classically the identity 
$\zeta(2)^2=2\zeta(2, 2)+4\zeta(3, 1)$. In our case, the same identity 
does not work, for example, if $p=2$, because the right side then is zero, 
but the left side is nonzero and is in fact transcendental. 

\subsection{Different identities with similar consequences!}
Classically, the multizeta values form a graded algebra.
More precisely, on account of the sum shuffle identities,
the product of multizeta values is a sum of multizeta values all
of  the same weight equal to the the sum of the weights of multizeta values in
the product. We were driven to introduce the degenerate multizeta
values into our picture in order to preserve this key classical piece of structure
in the function field setting.
But it seems (see \cite{T2} for details) that the nondegenerate multizeta
values still form an algebra, on account of relations which have an appearance
quite different from the shuffle type identities.
We only state here the simplest case of how the 
failures of shuffle identities mentioned above are salvaged by different identities: we have 
$\zeta(2)^2=\zeta(4)$, when $p=2$ (as can be seen from 3.1) or when $q=3$ 
(as can be seen from 3.2 by direct Bernoulli calculation or by direct calculation using generating functions). The sum shuffle above then fails for the
simple reason that $\zeta(2, 2)$ is non-zero as can be seen by a straight 
calculation of degrees.

\subsection{Sum-shuffle relations using degenerate multizeta values} 
We 
now make explicit the sum-shuffle relations already mentioned in
\ref{subsubsection:MotivationDegenerate}.
The following notation seems to be best suited for this.

\subsubsection{}
A {\em linear preorder} $\rho$ in a set $X$ 
is a relation satisfying the following 
conditions for all $x,y,z\in X$:
\begin{itemize}
\item  $x\rho y$ and $y\rho z$ implies $x\rho z$ (transitivity)
\item $x\rho y$ or $y \rho x$ (comparability) (It implies $x\rho x$ 
(reflexivity)).
\end{itemize}
Consider the following further axiom:
\begin{itemize}
\item $x\rho y$ and $y\rho x$ implies $x=y$ (antisymmetry)
\end{itemize}
A linear preorder which has also the antisymmetry property
is a  total ordering.  In general, for any linear preorder
$\rho$, the relation ``$x\rho y$ and $y\rho x$'' is an equivalence relation.
(In other words, a linear preorder is the same thing as an equivalence relation
along with a total order on the equivalence classes.)

\subsubsection{} Given a linear preorder $\rho$ in $\{1,\dots,r\}$,
we define the subset $(\NN)^r_\rho\subset (\NN)^r$
to be that whose members are the  $r$-tuples $(n_1,\dots,n_r)$
such that $i\rho j$ if and only if $n_i\leq n_j$
for all $i,j\in \{1,\dots,r\}$.

\subsubsection{} Given a linear preorder $\rho$ in $\{1,\dots,r\}$ 
and $s\in \ZZ_+^r$ put
$$\zeta_\rho(s)=\sum_{d\in (\NN)^r_\rho}S_d(s).
$$
These are exactly the multizeta values studied in this paper,
albeit expressed in a more elaborate notation in which the same multizeta value
may have several different presentations.
\subsubsection{} When we sum over all
linear preorders,  we get (analog of the  classical sum shuffle identity)
$$\prod_i \zeta(s_i)=\sum_\rho \zeta_\rho(s)$$
for all $s=(s_1, \cdots, s_r)\in \ZZ_+^r$.

\subsubsection{}\label{subsubsection:FFSumShuffle}
More generally, we now give 
the explicit identity showing that
$\FF_p$-linear combinations of the $\zeta_\rho(s)$ form an algebra.
For $\nu=0,1$ fix $s^{(\nu)}\in \ZZ_+^{r_\nu}$ and a linear preorder
$\rho_\nu$ on $\{1,\dots,r_\nu\}$.  Put $r=r_0+r_1$ and let $s\in
\ZZ_+^{r}$ be the concatenation of $s^{(0)}$ and $s^{(1)}$.  Then we
have
$$\zeta_{\rho_0}(s^{(0)})\zeta_{\rho_1}(s^{(1)})=\sum_\rho \zeta_\rho(s)$$
where the sum is extended over all linear pre-orders $\rho$ of $\{1,\dots,r\}$
such that $\rho$ restricts on $\{1,\dots,r_0\}$ to the linear pre-order $\rho_0$,
and on $\{r_0+1,\dots,r_1\}$ to the linear pre-order $\rho_1$ shifted by $r_0$.

\subsection{Relations from power sum relations involving digit conditions}

We refer to \cite{T2} for many combinatorially interesting  non-classical
relations between the non-degenerate multizeta values. Here we only
give some examples to illustrate  many more complicated relations that
exist once we allow degenerate multizetas. Let $z$ denote the `totally 
degenerate' (with all degree equalities) multizeta.  

The claims in the following examples follow by direct manipulations from 
$$
S_d(aq+b)=\frac{1}{\ell_d^{aq+b}}\left(1+\sum_{j=1}^a(-1)^j{b+j-1\choose j}
\frac{[d]^{jq}}{[1]^j}\right), \ \ \mbox{if} \ 0< a, b <q, $$
which follows easily from e.g., \cite[5.6.3]{T} or \cite{T2} or 2.4 above.

\subsubsection{}
By low $s$ relations above, $z(s_1, \cdots, s_r)=z(\sum s_i)$, if 
$s_i=k_ip^{n_i}$, $\sum s_i=kp^n$, with $k_i, k\leq q$. For example, 
$$z(p)=z(p-i, i)=z(p-j-k, j, k)=z(1, ...., 1)$$ 
and similarly we can change depth and replace parts in non-totally
degenerate multizetas at the places degeneration is allowed by such
method.

\subsubsection{}
Let  $p=2$. If  $0<b<q/3$, $b$ is of form $4k-1$, then 
$S_d(q+b)^3=S_d(3q+3b)$, whereas if $0<b<q/2$ and $b$ is of form 
$4k+1$, then $S_d(q+b)^3=S_d(3q+b)S_d(1)^{2b}$. Summing over $d$ gives 
relations between degenerate multizetas. 

\subsubsection{} 

Let $p\equiv 1 \mod 4$, choose $q, k$ and $0<b<q$ 
with $b^2\equiv -1 \mod p$ and $q>kp-2-2b>0$. Then 
\begin{eqnarray*}
&&S_d(q+b)^2+S_d(q+1)^2S_d(1)^{2b-2}\\
&&\;\;\;+S_d(q+(kp-2-2b))S_d(1)^{q-kp+2+4b}+(p-3)S_d(1)^{2q+2b}=0.
\end{eqnarray*}
Summing over $d$ we again get multizeta identities. 

As should be clear from these examples, we can find a wide variety 
of relations involving  various weights and depths.

\subsubsection{Motivic nature of relations and transcendence properties} 

Recent works have proved \cite{ABP, P} the analog of Grothendieck's
conjecture that the relations between periods should come from motivic
relations. Using this, together with the description \cite{AT} (see 2.4)
of Carlitz zeta values,
it has been proved in \cite{CY} that all the algebraic relations
between the Carlitz zeta values come from 3.1 and 3.2, and this is
still true \cite{CPY} even if you consider all with varying $q$ (in
the same characteristic) together.  Similarly, all algebraic relations 
between Carlitz zeta values and gamma values at proper fractions for 
the gamma function for $\FF_q[t]$ \cite[Chap. 4]{T} are known, by 
recent results of \cite{CPY2, CPTY}. (For history of earlier
transcendence results, see \cite[10.3, 10.5]{T}).

As for the higher depth multizeta values, some transcendence results in very
special circumstances were proved in \cite[10.5]{T}, using
manipulations of sums and techniques of \cite{AT, Y, Y97}.
By applying to the new identities in \cite{T2} the transcendence results
of \cite{P} on logarithms, more results were obtained.
It should be possible now to prove many
more transcendence and independence results
for multizeta values using the period interpretation given in this paper,
along with the general  transcendence results of \cite{P} that are now
available.

In contrast to the classical case, it is
a nice description of
the full set of identities in higher depths, not whether they would be
motivic, which is the hard part in our case! A good full description
of all identities is still being sought. See \cite{T2} for some progress. 

{\bf Acknowledgments}: The second  author thanks the Ellentuck Fund
and the von Neumann Fund for their support during his stay at the
Institute for Advanced Study, Princeton during Spring 2008, where this
paper was finalized. He is supported by NSA grant 98230-08-1-0049. 
We thank the referee for the suggestions to improve the exposition.

\end{document}